\documentclass[11 pt,oneside,reqno,a4paper]{amsart} 

\usepackage{amsfonts,amssymb,amscd,amsmath,enumerate,verbatim,calc} 
\usepackage{float}
\usepackage{amsthm}
\usepackage{extarrows}
\usepackage{mathrsfs}
\newtheorem{theorem}{Theorem}[section]
\newtheorem{definition}{Definition}[section]

\newtheorem{lemma}[theorem]{Lemma}
\newtheorem{corollary}[theorem]{Corollary}
\newtheorem{proposition}[theorem]{Proposition}
\newtheorem{remark}[theorem]{Remark}
\usepackage{mathtools}
\usepackage[all,cmtip]{xy}
\newtheorem{example}[theorem]{Example}

\numberwithin{equation}{section}

\usepackage{amsfonts}
\newcommand{\field}[1]{\mathbb{#1}}          

\newcommand{\N}{\field{N}}

\renewcommand{\ker}{\textnormal{ker}}

\begin{document}
	
		\title{Compatible multiplicative Lie algebra}
	\author{$^1$Amit Kumar and $^2$Sumit Kumar Upadhyay \\
		{$^{1,2}$Department of Applied Sciences,\\ Indian Institute of Information Technology Allahabad\\Prayagraj, U. P., India} }	

	\thanks{$^1$amitiiit007@gmail.com $^2$upadhyaysumit365@gmail.com}

	\begin{abstract}  
		In this paper, we introduce the concept of compatible multiplicative Lie algebras and develop the cohomology theory for compatible multiplicative Lie algebras. Our primary focus is on establishing the relationship between the compatible second cohomology group and central extensions of compatible multiplicative Lie algebras.
		\end{abstract}
	\maketitle 
	
	\section{Introduction}
The concept of compatible algebras appears in various areas of mathematics and mathematical physics, including the study of integrable systems, quantum groups, and deformations of algebraic structures.  In the late 1970s, I. M. Gelfand and I. Ya. Dorfman introduced compatible Hamiltonian operators: structures allowing two Poisson brackets whose linear combinations are also Poisson brackets \cite{GD}. During the 1980s and 1990s, many mathematicians further advanced the study of such structures, particularly within the theory of bi-Hamiltonian systems. These systems - central to integrable models - often involve a pair of compatible Lie or Poisson brackets. The theory of compatible Lie algebras was developed as a natural extension of studies in integrable systems and Poisson geometry. In this context, two or more Lie brackets are defined on the same vector space, and they satisfy specific compatibility relations ensuring that any linear combination of these brackets also forms a Lie algebra. This notion of compatibility has since been extended to a wide range of algebraic structures, including associative algebras, non-associative algebras, bialgebras, superalgebras, and many others. Recently, Ladra et al. \cite{MCA} gave a classification of nilpotent compatible Lie algebras.

A multiplicative Lie algebra (introduced by G. Ellis in \cite{GJ}) is a triple $ (G,\cdot,\star) $,  where $ (G,\cdot) $ is a group,  and $\star$ is a binary operation on $G$ such that for all $x, y, z  \in  G$, the following identities hold:
\begin{enumerate}
	\item $ x\star x=1 $, 
	\item $ x\star(yz)=(x\star y){^y(x\star z)} $, 
	\item $ (x y)\star z= {^x(y\star z)} (x\star z) $,
	\item $ ((x\star y)\star{^yz})((y\star z)\star{^zx})((z\star x)\star {^xy})=1 $,
	\item $ ^z(x\star y)=({^zx}\star {^zy})$, 
\end{enumerate}
where $^zx$ denotes the element $zxz^{-1}$. It is evident that any Lie algebra (and thus any Poisson algebra) is a multiplicative Lie algebra. Therefore, a multiplicative Lie algebra can be regarded as a generalization of groups, Lie algebras, and Poisson algebras. Thus, it is natural to extend the concept of compatibility to multiplicative Lie algebras.

In this paper, we first introduce the concept of a compatible multiplicative Lie algebra and then present some very nice illustrative examples. We then introduce the notions of nilpotency and solvability within the context of compatible multiplicative Lie algebras.

The problem of classifying finite groups naturally led to the development of extension theory and cohomology theory. This fundamental idea has since inspired the creation of cohomology theories for other algebraic structures. Over time, mathematicians have introduced numerous such theories, including the Chevalley-Eilenberg cohomology, the Hochschild cohomology, the cohomology of skew braces, the cohomology of Lie algebras, and multiplicative Lie algebras \cite{AGNM, RLS}, among others. 

In this study, we also develop a cohomology theory for compatible multiplicative Lie algebras using the framework of compatible central extensions. A key result is the establishment of a one-to-one correspondence between the compatible second cohomology group and the set of compatible central extensions. Furthermore, we establish a relationship between the compatible second cohomology group and the classical second cohomology groups.

\section{compatible multiplicative Lie algebra}
We start the section by defining a compatible multiplicative Lie algebra. The readers may refer to \cite{GJ, RLS, MS} for the basics of multiplicative Lie algebras.
 \begin{definition}\label{D1}
  Let $(G,\cdot)$ be a group with multiplicative Lie algebra structures $\star_1$ and  $\star_2$. We say that $(G,\cdot, \star_1, \star_2)$ is a compatible multiplicative Lie algebra if the following conditions hold:
  \begin{enumerate}
  \item $ G\star_1 G $ and $ G\star_2 G $ are abelian.
  \item $ (x\star_1 y)\cdot(z \star_2 w)=(z \star_2 w)\cdot(x\star_1 y)$  for all $x, y, z, w \in G$.
  \item $ ((x\star_1 y){\star_2} {^yz})\cdot((x\star_2 y){\star_1} {^yz})\cdot((y\star_1 z)\star_2{^zx})\cdot((y\star_2 z)\star_1{^zx})\cdot((z\star_1 x)\star_2 {^xy})\cdot((z\star_2 x)\star_1 {^xy})=1 $  for all $x, y, z \in G$.
\end{enumerate} 
 \end{definition}
The following theorem is a direct consequence of \cite[Theorem 3.5]{DAS}.
 \begin{theorem}\label{prod}
 Let $(G,\cdot)$ be a group with  multiplicative Lie algebra structures $\star_1$ and  $\star_2$ such that $ G\star_1 G $ and $ G\star_2 G $ are abelian and $ (x\star_1 y)\cdot(z \star_2 w)=(z \star_2 w)\cdot(x\star_1 y)$  for all $x, y, z, w \in G$. Then the binary operation $\star_{11} $ on $G,$ given by $x \star_{11} y = {(x\star_1 y)}\cdot{(x\star_2 y)} $, defines a multiplicative Lie algebra structure on $G$ if and only if $ ((x\star_1 y){\star_2} {^yz})\cdot((x\star_2 y){\star_1} {^yz})\cdot((y\star_1 z)\star_2{^zx})\cdot((y\star_2 z)\star_1{^zx})\cdot((z\star_1 x)\star_2 {^xy})\cdot((z\star_2 x)\star_1 {^xy})=1 $.
\end{theorem}

   \begin{remark}
  \begin{enumerate}
  \item Let $(G,\cdot)$ be an abelian group with multiplicative Lie algebra structures $\star_1$ and  $\star_2$. Then $(G,\cdot, \star_1, \star_2)$ is compatible if $ ((x\star_1 y){\star_2} {^yz})\cdot((x\star_2 y){\star_1} {^yz})\cdot((y\star_1 z)\star_2{^zx})\cdot((y\star_2 z)\star_1{^zx})\cdot((z\star_1 x)\star_2 {^xy})\cdot((z\star_2 x)\star_1 {^xy})=1 $  for all $x, y, z \in G$.
  \item  Let $(G,\cdot, \star_1, \star_2)$ be a compatible multiplicative Lie algebra. Then the binary operation $\star_{ij} $ on $G$, defined by $x \star_{ij} y = {(x\star_1 y)^i}\cdot{(x\star_2 y)^j} $ for all $i,j\in \N$, gives a multiplicative Lie algebra structure on $G$.
  \item If $\star_1 = \star_2$, then $(G,\cdot, \star_1, \star_2)$ is compatible if and only if $G\star_2 G = G\star_1 G$ is abelian.
  \item  Now onwards, we simply refer to $G$ or $(G,\star_1, \star_2)$ as a compatible multiplicative Lie algebra, rather than stating $(G,\cdot, \star_1, \star_2)$ is a compatible multiplicative Lie algebra.
  \item If $\star_1$ and  $\star_2$ both are trivial, then we say $G$ is trivially compatible.
  \item Suppose $\star_1$ is trivial and  $\star_2$ is a  non-trivial multiplicative Lie algebra structure. Then in Theorem \ref{prod}, $\star_{11} = \star_2$ is obviously multiplicative Lie algebra structure on $G$. So all the assumptions of the Definition \ref{D1} are redundant. In this case, we say $G$ is automatically compatible.
  \end{enumerate}
  
 \end{remark}
\begin{example}\label{trivia}
    \begin{enumerate}
        \item Let $G$ be a nilpotent group of class $2$ with multiplicative Lie algebra structures $\star_1$ and $\star_2$. If $\star_1 = \star_2$, then the conditions of Definition \ref{D1} are satisfied (see \cite{DAS}).
        \item Consider a nilpotent group $G$ of class $2$, endowed with two multiplicative Lie algebra structures $\star_1$ and $\star_2$ such that the multiplicative Lie algebra structures induced on $G/[G,G]$ are trivial. Then, the tuple $(G, \star_1, \star_2)$ forms a compatible multiplicative Lie algebra (see \cite{DAS}).
    \end{enumerate}
\end{example}
Now, we will give some more interesting examples of compatible multiplicative Lie algebras. 
\begin{example}
 Let $G$ be an abelian group with a multiplicative Lie algebra structure $\star$, and let $I$ be an ideal of $G$. 
 Then, the group $G\oplus I$ is a multiplicative Lie algebra with respect to the following three binary operations (see \cite{NAS}):
 \begin{align*}
    (g, a)\star_1 (h,b) &:= (g\star h,a\star b) \\
    (g, a)\star_2 (h,b) &:= \big(g\star h, (g\star b)(a\star h)\big)\\
    (g, a)\star_3 (h,b) &:= \big(g\star h, (g\star b)(a\star h)(a\star b)\big).
    \end{align*}
    It is easy to see that $(G\oplus I, \star_1, \star_2)$, $(G\oplus I, \star_1, \star_3)$ and $(G\oplus I, \star_2, \star_3)$ are compatible multiplicative Lie algebras.
 \end{example}

\begin{example}
  Let us consider the group $A_4 =\langle a, b \mid a^3 =b^2 =1, aba=ba^2b\rangle$ together with multiplicative Lie algebra structures
$\star_1$ and \(\star_2\) determined by
$a\star_1 b=[a,b]$, and
$a\star_2 b=b={^a[a,b]}$, see \cite{AD}. By \cite{MS}, we know that $A_4\star_1A_4=V_4=A_4\star_2A_4.$ Since \(V_4\) is abelian, both
$A_4\star_1A_4$ and $A_4\star_2A_4$ are abelian. Moreover, for all \(x,y,z,w\in A_4\), we have 
$(x\star_2y)(z\star_3w)=(z\star_3w)(x\star_2y)$.
It remains to verify condition (3) for compatible multiplicative Lie algebra. It is easy to see that
$
x\star_2 y={}^a[x,y]
$
for all $x,y\in A_4$.
Then
$
(x\star_1y)\star_2{}^yz
=
{}^a[x\star_1y,{}^yz]
=
{}^a[[x,y],{}^yz]
$
while
$
(x\star_2y)\star_1{}^yz
=
{}^a[x, y]\star_1{}^yz
=
[{}^a[x, y],{}^yz].
$
Put
$
u=[x,y]\in V_4,
$
and
$
q={}^yz\in A_4.
$
Since $V_4$ is a normal subgroup of $A_4$, $q$ belongs to one of the three cosets $V_4$, $V_4a$, or $V_4a^2$. Hence $q=va^i$ for some $v\in V_4$ and $i\in \{0,1,2\}$. This implies  $quq^{-1}=a^iu(a^i)^{-1}$. Thus, we have 
$$
{}^a[u,q]=[{}^au,q]
\qquad\text{for every }u\in V_4,\ q\in A_4.
$$
Therefore, we obtain
$(x\star_1y)\star_2{}^yz = (x\star_2y)\star_1{}^yz$. Similarly, 
$((y\star_1 z)\star_2{^zx})=((y\star_2 z)\star_1{^zx})~ \text{and}~ ((z\star_1 x)\star_2 {^xy})=((z\star_2 x)\star_1 {^xy})$. Now, since each non-identity element of $V_4$ is of order $2$, we have 
$$ ((x\star_1 y){\star_2} {^yz})\cdot((x\star_2 y){\star_1} {^yz})\cdot((y\star_1 z)\star_2{^zx})\cdot((y\star_2 z)\star_1{^zx})\cdot((z\star_1 x)\star_2 {^xy})\cdot((z\star_2 x)\star_1 {^xy})=1 $$  for all $x, y, z \in A_4$.
Hence $(A_4,\star_1,\star_2)$ is a compatible multiplicative Lie algebra.

\end{example}

A natural question arises at this point: Are any two multiplicative Lie algebra structures within the same group compatible $?$ The answer is, in general, negative.
\begin{example}
 Consider the symmetric group $S_n$ on $n$ symbols. Then it has only two multiplicative Lie algebra structures: a trivial one and another is defined by the commutator \cite{AD}. If we take $\star_1=\star_2=[,]$, for $n\geq 4$, $S_n$ fails to become compatible multiplicative Lie algebra as $[S_n,S_n]=A_n$. For another example, see \cite[Example 2.13]{MCA}.
\end{example}
 We now introduce the concepts of subalgebra, ideal, and homomorphism,  which lead to the notion of a quotient for a compatible multiplicative Lie algebra.

 \begin{definition}
 \begin{enumerate}
     \item A subalgebra $H$ of a compatible multiplicative Lie algebra $G$ is a subgroup of $G$ such that $(x\star_1 y)\in H $ and $(z\star_2 w)\in H$ for all $x,y,z,w \in H$.
     \item An ideal $I$ of a compatible multiplicative Lie algebra $G$ is a normal subgroup of $G$ such that $(x\star_1 y)\in I $ and $(z\star_2 w)\in I$ for all $x,z \in G$ and $y,w \in I$.
     \item  Let $(G, \star_1, \star_2)$ and $(\tilde G , \tilde\star_1, \tilde\star_2)$ be  compatible multiplicative Lie algebras. Then a group homomorphism $\phi: G \to \tilde G$ is called a compatible multiplicative Lie algebra homomorphism if
     \begin{align*}
         \phi(x\star_1 y) = \phi(x) \tilde\star_1 \phi(y) \\
         \phi(z\star_2 w) = \phi(z) \tilde\star_2 \phi(w)
     \end{align*}
     for all $x,y,z,w \in G$. Moreover, an invertible compatible multiplicative Lie algebra homomorphism is called a compatible multiplicative Lie algebra isomorphism.
 \end{enumerate}
\end{definition}
 \begin{remark}
 \begin{enumerate}
     \item A subalgebra  of a compatible multiplicative Lie algebra $(G, \star_1, \star_2)$ is also a subalgebra  of $(G, \star_{11})$.
     \item An ideal of a compatible multiplicative Lie algebra $(G, \star_1, \star_2)$ is also an ideal of $(G, \star_{11})$.
     \item The above definition implies that $\phi(x\star_{11} y) = \phi(x) \tilde\star_{11} \phi(y)$ for all $x,y \in G$, but the converse need not be true.
     \begin{example}
    Consider the quaternion group $Q_8 = \{\pm 1,\pm i, \pm j, \pm k\}$ with multiplicative Lie algebra structures $\star_1 = \star_2 = [,]$. Additionally, consider Klein's four group $V_4 = \{1, a, b, ab\}$ with multiplicative Lie algebra structures $ a \star_1 b = a \star_2 b = a$. Thus, both $Q_8$ and $V_4$  become compatible multiplicative Lie algebras in view of Example \ref{trivia}. Define $\phi: Q_8\to V_4$ by 
    \begin{align*}
\phi(\pm 1)=1,~~  \phi(\pm i) =a,~~  \phi(\pm j) =b,~~    \phi(\pm k)=ab.  
    \end{align*}
 Then  $\phi$ is a group homomorphism. It is easy to see that $\phi([x,y]^2)= (\phi(x)\star_1 \phi(y))^2 =1$  for all $x,y\in Q_8$. But $\phi([i,j])=\phi(-1)=1 $ and $ \phi(i)\star_1\phi(j) = a$ implies that $\phi([i,j])\neq \phi(i)\star_1\phi(j)$.
     \end{example}
     \item The above definition does not imply that $\phi$ is a multiplicative Lie algebra homomorphism between $(G, \star_1)$ and $(\tilde G,\tilde\star_2)$ or $(G,\star_2)$ and $(\tilde G, \tilde\star_1)$.
 \end{enumerate}
\end{remark}

 One can ask whether if $(G, \star_1)$ is isomorphic to $(\tilde G, \tilde\star_1)$ and $(G, \star_2)$ is isomorphic to $(\tilde G, \tilde\star_2)$, it necessarily follows that $(G, \star_1, \star_2)$ isomorphic to $(\tilde G , \tilde\star_1, \tilde\star_2)$. The answer is no. This is not always true.
 \begin{example}
   See \cite[Example 2.14]{MCA}.  
 \end{example}
  \begin{proposition}
         If $I$ is an ideal of a compatible multiplicative Lie algebra $G$, then the quotient $G/I$ forms a compatible multiplicative Lie algebra with respect to the following operations:
     \begin{align*}
        xI\star_1 yI := (x\star_1 y )I \\
        zI\star_2 wI := (z\star_2 w )I,
     \end{align*}
     where $x,y,z,w\in G$.
\end{proposition}
\begin{remark}
    As usual, there is a surjective compatible multiplicative Lie algebra homomorphism $\pi: G\to G/I$.
\end{remark}
 \begin{proposition}
    Let $\phi: G \to \tilde G$ be a compatible multiplicative Lie algebra homomorphism. Then the kernel of $\phi$ is an ideal of $G$, and the image of $\phi$ is a subalgebra of $\tilde G$. All isomorphism theorems follow in a routine manner.
 \end{proposition}
 We are now set to present various key tools that are vital for advancing the notion of nilpotency within compatible multiplicative Lie algebras. These tools also exemplify subalgebras and ideals.
 \subsection{Nilpotent compatible multiplicative Lie algebras}
  
\begin{definition}
   For a  compatible multiplicative Lie algebra $G$, we define 
   \begin{align*}
       \mathscr{LZ}(G) = \{x\in G \mid (x\star_1 y )=1=(x\star_2 z )~ \text{for all}~ y,z\in G\}\\
       \mathscr{MZ}(G) = \{x\in G \mid (x\star_1 y )=[x,y],(x\star_2 z )=[x,z]~ \text{for all}~ y,z\in G\}
   \end{align*}
  and refer to it as the compatible Lie center and the compatible multiplicative Lie center, respectively.
 \end{definition}

 \begin{remark}
 \begin{enumerate}
 \item $\mathscr{LZ}(G)$ is contained in $LZ(G) = \{x\in G \mid (x\star_{11} y )=1~ \text{for all}~ y\in G\}$.
 \item Let $G$ be a compatible multiplicative Lie algebra. Then $Z(G)$, $\mathscr{LZ}(G)$, $\mathscr{MZ}(G)$,  and $[G,G]$  are ideals of $G$, where $Z(G)$ and $[G,G]$ denote the group center and commutator subgroup of $G$, respectively.
 \end{enumerate}
     
 \end{remark}

 \begin{definition}
Let $G$ be a compatible multiplicative Lie algebra. Consider the set $S = \{(x\star_1 y ),(z\star_2 w )\mid x,y,z,w\in G\}$. Then, the subgroup of $G$ generated by the set $S$ is denoted by $G\star^{12}G$, and we call it a compatible commutator of $G$. 
 \end{definition}
 \begin{proposition}
   Let $G$ be a compatible multiplicative Lie algebra. Then the compatible commutator $G\star^{12}G$ is an ideal of $G$.
 \end{proposition}
 \begin{proof}
 It is enough to show that $a\star_1 ((x\star_1 y )(z\star_2 w ))$ and $a\star_2 ((x\star_1 y )(z\star_2 w ))$ belong to $G\star^{12}G$ for all $a,x,y,z,w\in G$.  
 \begin{align*}
 a\star_1 ((x\star_1 y )(z\star_2 w )) &= (a\star_1 (x\star_1 y )) ^{(x\star_1 y )} (a\star_1 (z\star_2 w ))\\&= (a\star_1 (x\star_1 y )) (a\star_1 (z\star_2 w )) \in G\star^{12}G.
 \end{align*}
 Similarly, \begin{align*}
 a\star_2 ((x\star_1 y )(z\star_2 w )) &= (a\star_2 (x\star_1 y )) ^{(x\star_1 y )} (a\star_2 (z\star_2 w ))\\&= (a\star_2 (x\star_1 y )) (a\star_2 (z\star_2 w )) \in G\star^{12}G.
 \end{align*}
 
 \end{proof}
 
\begin{definition}
     For a compatible multiplicative Lie algebra $G$, we define the lower central series of $G$ inductively as follows:
     \begin{align*}
         G = \Gamma_1(G)\unrhd \Gamma_2(G)\unrhd \Gamma_3(G)\unrhd \cdots \unrhd \Gamma_{n+1}(G)\unrhd \cdots,
     \end{align*}
    where $\Gamma_2(G)=G\star^{12}G$ and $\Gamma_{n+1}(G)=\Gamma_n(G)\star^{12} G$ for $n\geq 1$.

    Next, we define the upper central series of $G$ inductively as follows:
    \begin{align*}
        1=\mathscr{Z}_0(G)\unlhd \mathscr{Z}_1(G)\unlhd \mathscr{Z}_2(G)\unlhd\cdots \unlhd \mathscr{Z}_{n+1}(G)\unlhd\cdots,
    \end{align*}
    where $\mathscr{Z}_1(G)=\mathscr{LZ}(G)$ and for $n\geq 0$
    \begin{align*}
        \mathscr{Z}_{n+1}(G)=\{x\in G \mid (x\star_1 y )\in \mathscr{Z}_n(G), (x\star_2 z )\in \mathscr{Z}_n(G) ~ \text{for all}~ y,z\in G \},
    \end{align*}
    equivalently $\mathscr{Z}_{n+1}(G)/\mathscr{Z}_n(G) = \mathscr{LZ}(G/\mathscr{Z}_n(G))$.
    \end{definition}
\begin{lemma}\label{L1}
    Let $G$ be a compatible multiplicative Lie algebra. Then
    $$\Gamma_n(G)/\Gamma_{n+1}(G)\subseteq \mathscr{LZ}(G/\Gamma_{n+1}(G)).$$
\end{lemma}
\begin{theorem}
    Let $G$ be a compatible multiplicative Lie algebra. Then $\Gamma_{n+1}(G)=1$  if and only if $\mathscr{Z}_n(G) =G$.
\end{theorem}
\begin{proof} Suppose $\Gamma_{n+1}(G)=1$  with $\Gamma_n(G)\neq 1$. Now by induction, we claim that $\Gamma_{(n+2-i)}(G)\subseteq \mathscr{Z}_{i-1}(G)$ holds for all $i=1,2,\cdots ,{n+1}$.

By assumption, the claim is true for $i=1$. Suppose the statement is true for $i=m$, that is, $\Gamma_{(n+2-m)}(G)\subseteq \mathscr{Z}_{m-1}(G)$. Then by Lemma \ref{L1} and $\mathscr{LZ}(G/\Gamma_{(n+2-m)}(G)) \subseteq \mathscr{LZ}(G/\mathscr{Z}_{m-1}(G))$, we have $\Gamma_{(n+2-(m+1))}(G)\subseteq \mathscr{Z}_{m}(G)$. Thus, $\Gamma_{1}(G) = G\subseteq \mathscr{Z}_{n}(G)$,  that is $\mathscr{Z}_{n}(G) = G$.  Similarly, one can also prove the converse part.
\end{proof}
\begin{definition}
    A compatible multiplicative Lie algebra $G$ is called nilpotent of class $n$ if $\Gamma_{n+1}(G)=1$ with $\Gamma_n(G)\neq 1$  or equivalently $\mathscr{Z}_n(G) =G$ with $\mathscr{Z}_{n-1}(G)\neq G$.
\end{definition}
\begin{remark}
    We denote $(G,\cdot, \star_1)$ and $(G,\cdot, \star_2)$ by $G_1$ and $G_2$ respectively. Observe that if we have a nilpotent compatible multiplicative Lie algebra $(G, \star_1, \star_2)$ of class $n$,  both $G_1$ and $G_2$
    become nilpotent multiplicative Lie algebras of class $n$ according to \cite{FP}.
\end{remark}
The following corollary shows that the only nilpotent compatible multiplicative Lie algebra with trivial compatible Lie center is
the trivial compatible multiplicative Lie algebra.
\begin{corollary}
    Let $G$ be a nontrivial nilpotent compatible multiplicative Lie algebra. Then $\mathscr{LZ}(G) \neq 1$.
\end{corollary}
\begin{proof}
 If  $\mathscr{LZ}(G) = 1$, either $G=1$ or the upper central series never reaches $G$. Hence, the result follows.  
\end{proof}
A significant class of groups that arises in many areas of
research (Galois theory, for example) is the collection of solvable groups. Here, we aim to present an analogous concept of solvability in the context of compatible multiplicative Lie algebras.

Let $G$ be a compatible multiplicative Lie algebra. Define $\Gamma^n(G)$ inductively as follows: 
\begin{align*}
         G = \Gamma^1(G)\unrhd \Gamma^2(G)\unrhd \Gamma^3(G)\unrhd \cdots \unrhd \Gamma^n(G)\unrhd \cdots,
\end{align*}
where $\Gamma^2(G)=G\star^{12}G$ and $\Gamma^{n+1}(G)= \Gamma^n(G)\star^{12}\Gamma^n(G)$ for $n\geq 1$. This descending chain of ideals of $G$ is termed as derived series of $G$.
\begin{definition}
    A compatible multiplicative Lie algebra $G$ is called solvable if the derived series of $G$ terminates at $1$ after finitely many terms. The smallest $n$ such that $\Gamma^{n+1}(G)= 1$ is called the derived length of $G$.
\end{definition}
\begin{remark}
    If $G$ is a nilpotent compatible multiplicative Lie algebra then it is solvable. However, the reverse implication does not necessarily hold.
\end{remark}
\begin{example}
Consider Klein's four group $V_4 = \{1, a, b, ab\}$ with multiplicative Lie algebra structures $ a \star_1 b = a \star_2 b = a$. Then, by Example \ref{trivia}, $(V_4, \star_1, \star_2)$ becomes a compatible multiplicative Lie algebra. It is easy to see that $(V_4, \star_1, \star_2)$ is solvable but not nilpotent.
\end{example}

\subsection{Extensions and cohomology for compatible multiplicative Lie algebras}
A short exact sequence  
\begin{align*}
    \mathscr{E}: 1\longrightarrow  H \xlongrightarrow \alpha G \xlongrightarrow \beta K \longrightarrow 1
\end{align*}
of compatible multiplicative Lie algebras is called a compatible extension of $H$ by $K$. We may treat $\alpha $ as the inclusion map $i$ without any loss.
\begin{definition}
   An extension $\mathscr{E}$ is called a compatible central extension if
\begin{enumerate}
    \item $\alpha(H)\cong H\subseteq Z(G)$, that is,  $\mathscr{E}$ is a central extension of groups.
    \item $\alpha(H)\cong H\subseteq \mathscr{LZ}(G)  $.
\end{enumerate}   
In other words, a compatible extension $\mathscr{E}$ is called a compatible central extension if $ H\subseteq Z(G)\cap \mathscr{LZ}(G)  $.
\end{definition}
Let
\begin{align*}
    \mathscr{E'}: 1\longrightarrow  H \xlongrightarrow {i} G' \xlongrightarrow {\beta'} K \longrightarrow 1
\end{align*}
be another compatible central extension. We say that $\mathscr{E}$ and $\mathscr{E'}$ are equivalent if there exists a compatible multiplicative Lie algebra isomorphism  $\gamma: G\to G'$ such that the following diagram is commutative 
\[
	\xymatrix{
	 1\ar[r] &	H\ar[d]_{I_H}\ar[r]^{i~} & G \ar[r]^{\beta~} \ar[d]_{\gamma}  & K  \ar[d]_{I_K} \ar[r] & 1 \\		1\ar[r] & H\ar[r]^{i~}& G' \ar[r]^{\beta'~} & K \ar[r] & 1 .}
	\]

Let $t$ be a section of the compatible central extension $\mathscr{E}$. Then we have a map $f^t:K\times K \to H$ such that
\begin{align}
    &f^t(1,x)=f^t(x,1)=1 \\
    &f^t(x,y)f^t(xy,z)=f^t(y,z)f^t(x,yz)
\end{align}
for all $x,y,z\in K$. We know that every element of $G$ is uniquely written as $at(x)$ for some $a\in H$ and $x\in K$. The group operation on $G$ is given by 
\begin{align}
    (at(x))(bt(y)) = abf^t(x,y)t(xy).
\end{align}
The compatible multiplicative Lie algebra structure on $G$ is given by
\begin{align}
    (at(x))\star_i(bt(y)) = h_i^t(x,y)t(x\star_iy)
\end{align}
such that the maps $h^t_i:K\times K \to H$, $i=1,2$, satisfy  the following equations: 
\begin{equation}\label{H1}
    h_i^t(x,1)=h_i^t(1,x)=h_i^t(x,x)=1 
\end{equation}
\begin{equation}\label{H2}
    h_i^t(x,yz)=h_i^t(x,y)h_i^t(x,z)f^t(y^{-1},y)^{-1}f^t(y,x\star_i z) f^t(y(x\star_i z),y^{-1})f^t(x\star_i y,{^y(x\star_i z)}) 
\end{equation}
\begin{equation}\label{H3}
    h_i^t(xy,z)=h_i^t(y,z)h_i^t(x,z)f^t(x^{-1},x)^{-1}f^t(x,y\star_i z) f^t(x(y\star_i z),x^{-1})f^t({^x(y\star_i z)},(x\star_i z))
\end{equation}
\begin{equation}\label{H4}
\hspace{-1cm} h_i^t(y\star_i x,{^xz})h_i^t(x\star_i z,{^zy})h_i^t(z\star_i y,{^yx})f^t((y\star_i x)\star_i {^xz},(x\star_i z)\star_i {^zy}) f^t(((y\star_i x)\star_i {^xz})((x\star_i z)\star_i {^zy}),(z\star_i y)\star_i {^yx})=1
\end{equation}
\begin{equation}\label{H5}
    h_i^t({^zx},{^zy})=h_i^t(x,y)f^t(z,x\star y)f^t(z^{-1},z)^{-1}f^t(z(x\star_i y),z^{-1})
\end{equation}
for all $x,y,z\in K$. To see more details, we refer to \cite{RLS}. 
Further, using the Definition \ref{D1}, the following equations also hold:
\begin{align*}
    &((at(x))\star_i(bt(y)))((ct(z))\star_i(dt(w)))
    =  ((ct(z))\star_i(dt(w)))((at(x))\star_i(bt(y)))((ct(z)).\\
\implies &(h_i^t(x,y)t(x\star_iy))(h_i^t(z,w)t(z\star_iw))=(h_i^t(z,w)t(z\star_iw))(h_i^t(x,y)t(x\star_iy)).\\
\implies &t(x\star_iy)t(z\star_iw)=t(z\star_iw)t(x\star_iy).\\
\implies &f^t(x\star_iy,z\star_iw)t((x\star_iy)(z\star_iw))=f^t(z\star_iw,x\star_iy)t((z\star_iw)(x\star_iy)).
\end{align*}
\vspace{-1cm} 
\begin{align}\label{Ab1}
\hspace{-7.5cm}\implies f^t(x\star_iy,z\star_iw)=f^t(z\star_iw,x\star_iy).
\end{align}
Similarly, we have
\begin{align*}
    &((at(x))\star_1(bt(y)))((ct(z))\star_2(dt(w)))
    =  ((ct(z))\star_2(dt(w)))((at(x))\star_1(bt(y)))((ct(z)).
 \end{align*}
\vspace{-1cm} 
\begin{align}\label{Ab2}
\hspace{-7.5cm}\implies f^t(x\star_1y,z\star_2w)=f^t(z\star_2w,x\star_1y).
\end{align}
Next,
\begin{align*}
 & h^t_2(x\star_1y,{^yz})t((x\star_1 y)\star_2 {^yz}) h^t_1(x\star_2y,{^yz})t((x\star_2 y)\star_1 {^yz})h^t_2(y\star_1z,{^zx})t((y\star_1 z)\star_2 {^zx}) \\
 & h^t_1(y\star_2z,{^zx})t((y\star_2 z)\star_1 {^zx})h^t_2(z\star_1x,{^xy})t((z\star_1 x)\star_2 {^xy}) h^t_1(z\star_2x,{^xy})t((z\star_2 x)\star_1 {^xy}) =1.\\
 \implies & h^t_2(x\star_1y,{^yz})h^t_1(x\star_2y,{^yz})h^t_2(y\star_1z,{^zx})h^t_1(y\star_2z,{^zx})h^t_2(z\star_1x,{^xy})h^t_1(z\star_2x,{^xy}) \\
 & t((x\star_1 y)\star_2 {^yz}) t((x\star_2 y)\star_1 {^yz})t((y\star_1 z)\star_2 {^zx})t((y\star_2 z)\star_1 {^zx})t((z\star_1 x)\star_2 {^xy}) t((z\star_2 x)\star_1 {^xy}) =1.\\
 \implies & h^t_2(x\star_1y,{^yz})h^t_1(x\star_2y,{^yz})h^t_2(y\star_1z,{^zx})h^t_1(y\star_2z,{^zx})h^t_2(z\star_1x,{^xy})h^t_1(z\star_2x,{^xy}) \\
 & f^t((x\star_1 y)\star_2 {^yz},(x\star_2 y)\star_1 {^yz})f^t(((x\star_1 y)\star_2 {^yz})((x\star_2 y)\star_1 {^yz}),(y\star_1 z)\star_2 {^zx}) \\
 &f^t(((x\star_1 y)\star_2 {^yz})((x\star_2 y)\star_1 {^yz})((y\star_1 z)\star_2 {^zx}),(y\star_2 z)\star_1 {^zx})\\
 &f^t(((x\star_1 y)\star_2 {^yz})((x\star_2 y)\star_1 {^yz})((y\star_1 z)\star_2 {^zx})((y\star_2 z)\star_1 {^zx}),(z\star_1 x)\star_2 {^xy})\\
 &f^t(((x\star_1 y)\star_2 {^yz})((x\star_2 y)\star_1 {^yz})((y\star_1 z)\star_2 {^zx})((y\star_2 z)\star_1 {^zx})((z\star_1 x)\star_2 {^xy}),(z\star_2 x)\star_1 {^xy})\\
 &t(((x\star_1 y){\star_2} {^yz})((x\star_2 y){\star_1} {^yz})((y\star_1 z)\star_2{^zx})((y\star_2 z)\star_1{^zx})((z\star_1 x)\star_2 {^xy})((z\star_2 x)\star_1 {^xy})) = 1.\\
  \implies & h^t_2(x\star_1y,{^yz})h^t_1(x\star_2y,{^yz})h^t_2(y\star_1z,{^zx})h^t_1(y\star_2z,{^zx})h^t_2(z\star_1x,{^xy})h^t_1(z\star_2x,{^xy}) \\
 & f^t((x\star_1 y)\star_2 {^yz},(x\star_2 y)\star_1 {^yz})f^t(((x\star_1 y)\star_2 {^yz})((x\star_2 y)\star_1 {^yz}),(y\star_1 z)\star_2 {^zx}) \\
 &f^t(((x\star_1 y)\star_2 {^yz})((x\star_2 y)\star_1 {^yz})((y\star_1 z)\star_2 {^zx}),(y\star_2 z)\star_1 {^zx})\\
 &f^t(((x\star_1 y)\star_2 {^yz})((x\star_2 y)\star_1 {^yz})((y\star_1 z)\star_2 {^zx})((y\star_2 z)\star_1 {^zx}),(z\star_1 x)\star_2 {^xy})\\
 &f^t(((x\star_1 y)\star_2 {^yz})((x\star_2 y)\star_1 {^yz})((y\star_1 z)\star_2 {^zx})((y\star_2 z)\star_1 {^zx})((z\star_1 x)\star_2 {^xy}),(z\star_2 x)\star_1 {^xy})=1.\\
 \implies & h^t_2(x\star_1y,{^yz})h^t_1(x\star_2y,{^yz})h^t_2(y\star_1z,{^zx})h^t_1(y\star_2z,{^zx})h^t_2(z\star_1x,{^xy})h^t_1(z\star_2x,{^xy}) \\
 & f^t((x\star_1 y)\star_2 {^yz},(x\star_2 y)\star_1 {^yz})f^t(((x\star_1 y)\star_2 {^yz})((x\star_2 y)\star_1 {^yz}),(y\star_1 z)\star_2 {^zx}) \\
 &f^t(((x\star_1 y)\star_2 {^yz})((x\star_2 y)\star_1 {^yz})((y\star_1 z)\star_2 {^zx}),(y\star_2 z)\star_1 {^zx})
\end{align*}
 \vspace{-1cm}
\begin{align}\label{MixJ}
\hspace{.8cm}f^t(((z\star_2 x)\star_1 {^xy})^{-1}((z\star_1 x)\star_2 {^xy})^{-1},(z\star_1 x)\star_2 {^xy})f^t(((z\star_2 x)\star_1 {^xy})^{-1},(z\star_2 x)\star_1 {^xy})=1.
\end{align}
\begin{remark}
A compatible central extension $\mathscr{E}$ with a section $t$ gives us a triple $(f^t,h^t_1,h^t_2)$ that satisfies the aforementioned equations.  
\end{remark}
\begin{definition}
    A compatible $2-$cocycle of a compatible multiplicative Lie algebra $K$ with coefficients in an abelian group $H$ which is  trivially compatible is a triple $(f,h_1,h_2)$, where $f\in Z^2(K, H)$ and $h_1,h_2:K\times K\to H$ satisfy the equations $(\ref{H1}),(\ref{H2}),(\ref{H3}),(\ref{H4}),(\ref{H5}),(\ref{Ab1}),(\ref{Ab2})$, and $(\ref{MixJ})$.
\end{definition}

Conversely, suppose $(f,h_1,h_2)$ is a compatible $2-$cocycle of a compatible multiplicative Lie algebra $K$ with coefficients in an abelian group $H$ which is  trivially compatible. Take $G=H\times K$. We define a group operation on $G$ as 
\begin{align*}
    (a,x)(b,y)=(abf(x,y),xy)
\end{align*}
and two binary operations on $G$ by
\begin{align*}
   &(a,x)\tilde\star_1(b,y)=(h_1(x,y),x\star_1y) \\
   &(a,x)\tilde\star_2(b,y)=(h_2(x,y),x\star_2y).
\end{align*}
It can be easily checked that $\tilde\star_1$ and $\tilde\star_2$ are multiplicative Lie algebra structures on $G$ using the equations $(\ref{H1}),(\ref{H2}),(\ref{H3}),(\ref{H4}),(\ref{H5})$. By the use of equations (\ref{Ab1}), (\ref{Ab2}) and (\ref{MixJ}), it is easy to see that $(G,\tilde\star_1,\tilde\star_2)$ is a compatible multiplicative Lie algebra such that 
\begin{align*}
    \mathscr{N}: 1\longrightarrow  H \xlongrightarrow i G \xlongrightarrow p K \longrightarrow 1
\end{align*}	
is a compatible central extension of $H$ by $K$, where $i$ is the first inclusion and $p$ is the second projection map. The map $t:K\to G$ given by $t(x)=(1,x)$ is a section such that $f^t=f$, $h^t_1=h_1$ and 
$h^t_2=h_2$. We call $\mathscr{N}$ as the compatible central extension
associated to the compatible $2-$cocycle $(f,h_1,h_2)$.

If we have another section $s$ of $\mathscr{E}$. Then there is a map $g:K\to H$ with $g(1)=1$ such that $s(x)=g(x)t(x)$ for each $x\in K$ and 
\begin{align*}
    &f^s(x,y)=g(x)g(xy)^{-1}g(y)f^t(x,y)\\
    &h^s_1(x,y)=g(x\star_1y)^{-1}h^t_1(x,y)\\
    &h^s_2(x,y)=g(x\star_2y)^{-1}h^t_2(x,y)
\end{align*}
for all $x,y\in K$. This motivates the following definition:
\begin{definition}
    The compatible $2-$cocycles $(f,h_1,h_2)$ and $(f',h'_1,h'_2)$ are said to be equivalent if there exists an identity preserving map $g: K\to H$ such that 
    \begin{align*}
    &f'(x,y)=g(x)g(xy)^{-1}g(y)f(x,y)\\
    &h'_1(x,y)=g(x\star_1y)^{-1}h_1(x,y)\\
    &h'_2(x,y)=g(x\star_2y)^{-1}h_2(x,y)
\end{align*}
for all $x,y\in K$.
\end{definition}
It is easy to see that the set $Z^2_C(K, H)$ of all compatible $2-$cocycles of $K$ with coefficient in $H$ is an abelian group with respect to the following operation:
\begin{align*}
    (f,h_1,h_2)(f',h'_1,h'_2)=(ff',h_1h'_1,h_2h'_2).
\end{align*}
If $g:K\to H$ is a given identity preserving map, then the triple $(\delta g,g^{\star_1},g^{\star_2})$ of maps from $K\times K$ to $H$ given by 
\begin{align*}
   \delta g((x,y))&= g(x)g(xy)^{-1}g(y)\\
   g^{\star_1}((x,y))&=g(x\star_1y)^{-1}\\
   g^{\star_2}((x,y))&=g(x\star_2y)^{-1}
\end{align*}
can be easily seen as a member of $Z^2_C(K, H)$.

Let $\mathcal{F}(K, H)$ denote the group of identity-preserving maps from $K$ to $H$. Then we have a homomorphism $\lambda:\mathcal{F}(K,H)\to Z^2_C(K,H)$ given by $\lambda(g)=(\delta g,g^{\star_1},g^{\star_2})$. The image of $\lambda$ is called the group of compatible $2-$coboundaries of $K$ with coefficients in $H$, denoted by $B^2_C(K, H)$. The quotient group $Z^2_C(K, H)/B^2_C(K, H)$ is called the compatible second cohomology of $K$ with coefficients in $H$ and is denoted by $H^2_C(K, H)$.
Thus, we have the following exact sequence of abelian groups:
\begin{align}
  1\longrightarrow Hom(K,H) \xlongrightarrow i \mathcal{F}(K,H) \xlongrightarrow \lambda Z^2_C(K,H) \xlongrightarrow \nu H^2_C(K,H) \longrightarrow 1,  
\end{align}
where $\nu$ is the quotient map. Also, we have a homomorphism $p:Z^2_C(K,H)\to Z^2(K,H)$ given by $p((f,h_1,h_2))=f$. It is clear that $p(B^2_C(K,H))\subseteq B^2(K,H)$. Therefore $p$ induces a homomorphism $\Tilde{p}:H^2_C(K,H) \to H^2(K,H)$ given by $\Tilde{p}\big(\overline{(f,h_1,h_2)}\big)=\overline{f}$.

We are now ready to tie together all equivalence classes of compatible central extensions of $H$ by $K$ and the elements of the group $H^2_C(K,H)$.
\begin{theorem}
    There is a bijective correspondence between the equivalence classes of compatible central extensions of $H$ by $K$ and the compatible second cohomology group $H^2_C(K,H)$. 
\end{theorem}
\begin{proposition}\label{kernel}
The kernel of the homomorphism $\Tilde{p}:H^2_C(K,H) \to H^2(K,H)$ is isomorphic to the subgroup
\begin{align*}
 Hom\times Hom=\{(\phi_1,\phi_2)\in Hom(\wedge^2K/J_1,H)\times Hom(\wedge^2K/J_2,H)\mid \phi_1(X_2)\phi_2(X_1)=1 \}   
\end{align*} 
of $Hom(\wedge^2K/J_1,H)\times Hom(\wedge^2K/J_2,H)$, where  $\wedge^2K$ is the non abelian exterior square of $K$, $J_1$ is the subgroup of $\wedge^2K$ generated by the set 
\begin{align*}
\{((x\star_1 y)\wedge{^yz})((y\star_1 z)\wedge{^zx})((z\star_1 x)\wedge {^xy})\mid x,y,z \in K \},    
\end{align*}
$J_2$ is the subgroup of $\wedge^2K$ generated by the set 
\begin{align*}
  \{((x\star_2 y)\wedge{^yz})((y\star_2 z)\wedge{^zx})((z\star_2 x)\wedge {^xy})\mid x,y,z \in K \},  
\end{align*}
$X_2 = ((x\star_2 y)\wedge{^yz})((y\star_2 z)\wedge{^zx})((z\star_2 x)\wedge {^xy}) $  and 
$X_1 = ((x\star_1 y)\wedge{^yz})((y\star_1 z)\wedge{^zx})((z\star_1 x)\wedge {^xy}) $ for all $x,y,z \in K$.
\end{proposition}
\begin{proof}
 Every element of $\ker(\Tilde{p})$ is uniquely expressible as $(1,h_1,h_2)+B^2_C(K,H)$, where $h_1,h_2:K\times K \to H$ are the maps satisfying the following equations for $i=1,2$ and $x,y,z\in K$ :
 \begin{align*}
    &h_i(x,1)=h_i(1,x)=h_i(x,x)=1.\\
    &h_i(x,yz)=h_i(x,y)h_i(x,z).\\
    &h_i(xy,z)=h_i(y,z)h_i(x,z).\\
    &h_i({^zx},{^zy})=h_i(x,y).\\
    &h_i(x\star_i y,{^yz})h_i(y\star_i z,{^zx})h_i(z\star_i x,{^xy})=1.\\
    h_1(x\star_2y,{^yz})&h_1(y\star_2z,{^zx})h_1(z\star_2x,{^xy})h_2(x\star_1y,{^yz})h_2(y\star_1z,{^zx})h_2(z\star_1x,{^xy})=1.
\end{align*}
This induces homomorphisms $\tilde h_1:\wedge^2K/J_1\to H$ and $\tilde h_2:\wedge^2K/J_2\to H$ such that 
\begin{align*}
 \tilde h_1(((x\star_2 y)\wedge{^yz})((y\star_2 z)\wedge{^zx})((z\star_2 x)\wedge {^xy}))\tilde h_2(((x\star_1 y)\wedge{^yz})((y\star_1 z)\wedge{^zx})((z\star_1 x)\wedge {^xy}))=1   
\end{align*}
for all $x,y,z \in K$. Thus, $\ker(\Tilde{p})$ is isomorphic to the group $Hom\times Hom$.
\end{proof}
Consider a non-abelian group $K$ for which the commutator subgroup $[K, K]$ is abelian. Then, the structure $(K,\star_1,\star_2)=(K,[,],[,])$ forms a compatible multiplicative Lie algebra. Such a structure is termed a compatible improper multiplicative Lie algebra.
\begin{lemma}
 Let $K$ be a compatible improper multiplicative Lie algebra. Then the  homomorphism $\Tilde{p}:H^2_C(K,H) \to H^2(K,H)$ is surjective.
\end{lemma}
\begin{proof}
    Assume $f+B^2(K,H)\in  H^2(K,H)$. Then, there exists a central extension of groups
    \begin{align*}
 \mathscr{C}: 1\longrightarrow  H \xlongrightarrow i G \xlongrightarrow \mu K \longrightarrow 1
    \end{align*}
    such that the corresponding $2-$cocycle is $f$. Define the compatible  multiplicative Lie algebra on $G$ with respect to the section $t$ as follows:
    \begin{align*}
(at(x))\star_i(bt(y))=f(x,y)f(x^{-1},x)^{-1}f(xy,x^{-1})f(y^{-1},y)^{-1}f(xyx^{-1},y^{-1})t([x,y]) 
\end{align*}
Thus, the central extension $\mathscr{C}$ is also a compatible central extension of $H$ by the compatible improper multiplicative Lie algebra $K$. This implies that 
 \begin{align*}   
     (at(x))\star_i(bt(y)) = h_i^t(x,y)t([x,y])  ~~ \text{for}~~ i=1,2.
 \end{align*}
 Hence we get
  \begin{align*}   
    h_1^t(x,y) = f(x,y)f(x^{-1},x)^{-1}f(xy,x^{-1})f(y^{-1},y)^{-1}f(xyx^{-1},y^{-1}) =  h_2^t(x,y).
 \end{align*}
Thus, given a $2-$cocycle  $f$, we get a compatible $2-$cocycle $(f,h_1^t,h_2^t)$ such that  $\Tilde{p}((f,h_1^t,h_2^t))=f$. Hence $\Tilde{p}$ is surjective.
\end{proof}
\begin{corollary}
Let $K$ be a group with a trivial Schur multiplier. Then $H^2_C(K, H)$ is isomorphic to $H^2(K, H)\times (Hom\times Hom)^{\times}$, where $K$ is a compatible improper multiplicative Lie algebra and
\begin{align*}
 (Hom\times Hom)^{\times} &=  \{(\phi_1,\phi_2)\in Hom(\wedge^2K,H)\times Hom(\wedge^2K,H)\mid (\phi_1\phi_2)(([x, y]\wedge{^yz})\\ &([y, z]\wedge{^zx})([z,x]\wedge {^xy}))=1 ~\text{for all}~ x,y,z\in K \}.  
\end{align*}
\end{corollary}
\begin{proof}
    The Schur multiplier of $K$ is trivial, implying that the subgroup $J_1=J_2$ of $\wedge^2K$ is trivial \cite[Corollary 3.5]{RLS}. Hence, by Proposition \ref{kernel}, we obtain the following short exact sequence:
    \begin{align*}
     1\longrightarrow   (Hom\times Hom)^{\times} \xlongrightarrow i H^2_C(K,H) \xlongrightarrow {\tilde{p}} H^2(K,H) \longrightarrow 1. 
     \end{align*}
     The map $t: H^2(K,H) \to H^2_C(K,H)$ given by $t(\overline{f})=\overline{(f,h_1,h_2)}$ is a splitting, where 
     \begin{align*}   
    h_1(x,y) = f(x,y)f(x^{-1},x)^{-1}f(xy,x^{-1})f(y^{-1},y)^{-1}f(xyx^{-1},y^{-1}) =  h_2(x,y).
 \end{align*}
 Hence the result follows.
\end{proof}

\noindent{\bf Acknowledgement:} 
The first named author sincerely thanks IIIT Allahabad and the University Grant Commission (UGC), Govt. of India, New Delhi, for a research fellowship. The second named author is thankful to the National Board for Higher Mathematics (NBHM) for the financial support for the project ``Linear Representation of Multiplicative Lie Algebra" (02011/19/2023/NBHM (R.P)/R~\& D-II/5954).

We also thankful to the anonymous reviewer for the valuable questions and constructive comments about the existence of self-compatible multiplicative Lie algebras beyond those arising from abelian groups and nilpotent groups of class $2$, and the possibility of a universal compatible central extension. These questions have helped us clarify the scope of our results and identify directions for further investigation.

\noindent \textbf{Statements and Declarations.}

\smallskip

\noindent {\bf Author contributions.} The authors contributed equally to the preparation of this article.\\
\noindent  {\bf Funding.} The authors confirm that they did not receive any funding, grants, or support during the preparation of this manuscript.\\
\noindent  {\bf Data availability.} Data sharing does not apply to this article, as no datasets were analyzed during its preparation.\\
\noindent   {\bf Conflict of interests.} The authors affirm that they have no financial interests to disclose.


\begin{thebibliography}{9}

\bibitem{AGNM}
 A. Bak, G. Donadze, N. Inassaridze, M. Ladra,
Homology of multiplicative Lie rings, J. Pure Appl. Algebra, 208 (2007), 761-777.

\bibitem{GJ}
G. J. Ellis, On five well-known commutator identities, J. Aust. Math. Soc. (Series A), 54 (1993), 1-19.



\bibitem{GD} 
 I. M. Gelfand, and  I. Ya. Dorfman, Hamiltonian operators and algebraic structures related to them. Funct. Anal. Appl., 13 (4) (1979), 248–262.
	
\bibitem{AD} 
A. Kumar,  D. Pal, S. Kushwaha, and S. K. Upadhyay, The Schur multiplier of some finite Multiplicative Lie algebras, J. Lie Theory, 34 (2) (2024), 423-436.

\bibitem{MCA}
 M. Ladra, B. L. d. Cunha, and S.A. Lopes, A classification of nilpotent compatible Lie algebras. Rend. Circ. Mat. Palermo, II. Ser 74, 70 (2025). 

\bibitem{RLS}
 R. Lal and S. K. Upadhyay, Multiplicative Lie algebra and Schur multiplier, J. Pure Appl. Algebra, 223 (9) (2019), 3695-3721.



\bibitem{NAS} 
 N. K. Maurya, A. Kumar, and S. K. Upadhyay, Excision and idealization of a multiplicative Lie algebra, Ricerche mat, 75 (4) (2026), 2419–2434. 
	  
\bibitem{DAS}
 D. Pal,  A. Kumar, and  S. K. Upadhyay, Multiplicative Lie algebra structure on nilpotent group of class $2$, Georgian Math. J., 33  (2) (2026),  301-306. 
	

	
\bibitem{MS}
	M. S. Pandey and S. K. Upadhyay, Classification of multiplicative Lie algebra structures on a finite group, Colloq. Math., 168 (2022), 25-34.
	
\bibitem{FP}
	F. Point and P. Wantiez, Nilpotency criteria for multiplicative Lie algebra, J. Pure Appl. Algebra, 111 (1996), 229-243.
	




\end{thebibliography}
\end{document}